**BAIBEKOV S.N., DOSSAYEVA A.A.**

*R&D center "Kazakhstan engineering"*
*Astana c., Kazakhstan*

# DEVELOPMENT OF THE MATRIX OF PRIMES AND PROOF OF AN INFINITE NUMBER OF PRIMES -TWINS

---

**Abstract:** *This paper is devoted to the theory of prime numbers. In this paper we first introduce the notion of a matrix of prime numbers. Then, in order to investigate the density of prime numbers in separate rows of the matrix under consideration, we propose a number of lemmas and theorems that, together with the Dirichlet and Euler theorems, make it possible to prove the infinity of prime twins.*

**Keywords:** *Prime numbers, prime-twins, composite numbers, integers, algorithm, arithmetical progression, matrix of prime numbers, special factorial.*

---

## 1. Introduction

In 1912, during 5-th International Mathematical Congress, which was held at the University of Cambridge, it was noted that there is a number of open problems in the field of prime numbers. One of them is a problem of proving the infinity of number of prime twin numbers. This problem, which has not been solved for more than 2000 years as of today, is now referred as the second Landau problem.

The aim of this work is to prove the infinity of prime twins.

In [1] and [2] we proposed two versions of the analytic proof of the infinity of prime twin numbers, which, we believe, give a correct result about the infinity of twins. In this paper we propose another new proof of the infinity of twins, which, if possible, turned out to be more "elementary" and more explicit. We believe that the successful use of the fundamental works of world famous scientists (Eratosthenes, Euler, Dirichlet and others) allowed us to do this.

First, for convenience, as in [1] and [2], we introduce the following notation. As known, a consecutive multiplication of natural numbers is called a factorial:

$\prod_{i=1}^{n} i = n!$ In the future, sequential multiplication of primes will occur quite often, so for such cases we use the notation:

$$2 * 3 * 5 * 7 * 11 * \ldots * \ p_n = \prod_{i=1}^{n} p_i = p_n!'.$$

Here $p_i$ is a prime number with a serial number $i$. A combination of symbols $p_n!'$ represent consistent multiplication of prime numbers from 2 to $p_n$. It will be called a **special factorial (or primorial) of prime number $p_n$.** For example, $p_4!'$ - is a special factorial of fourth prime number $p_4 = 7$ or $p_4!' = 7!' = 2 * 3 * 5 * 7 = 210$.

To achieve the goal set for the proposed work, we develop matrices of prime numbers.

## 2. Matrix of prime numbers

Let us represent the set of natural numbers (except for number 1) in a form of a matrix family $A_k$ with elements of $a(k, i_k, j_k)$, where $i_k$ is a serial number of rows, $j_k$ is a serial number of columns and $k$ – is a serial number of matrix $A_k$.

$$a(k, i_k, j_k) = (i_k + 1) + p_k!'(j_k - 1) = l_k + D_k(j_k - 1), \qquad (1)$$

where $j_k = 1, 2, \ldots, \infty$; $i_k = 1, 2, \ldots, p_k!'$. $l_k = i_k + 1$ - is the first number of plurality of numbers that are in $i_k$ row. Note, index $k$ of the parameters $i$ and $j$ denotes affiliation of these parameters to matrix $A_k$.

It appears from the equation (1) that the sequence of numbers in any row of the matrix $A_k$, is an arithmetic progression with first term equal to $l_k = i_k + 1$. The difference in this progression is $D_k = p_k!'$.

In this case, the maximum number of rows of matrix $A_k$ should be equal to a special factorial $p_k!'$, i.e.:

$$i_{k,max} = p_k!'$$

The number of columns can be arbitrarily large up to infinity (Fig. 1).

As known, the number 1 is not a prime number, otherwise all numbers that are multiples of 1 would be composite numbers. The number 1 is also not a composite number, since it is not divisible by other numbers. That is why the number 1 is located separately in the upper left corner of the matrices. In Fig. 1, only four matrices $A_1$, $A_2$, $A_3$ and $A_4$ are shown for clarity. The number of rows of other matrices exceeds the paper size, so they are not shown in Fig 1.

It follows from (1) that each matrix from the family of matrices $A_k$, where $k = 1, 2, 3, \ldots, \infty$, differs from each other, primarily by the number of rows and the difference in arithmetic progression. For example, in case of matrix $A_{k-1}$, the arithmetic progression constant will be equal to $D_{k-1} = p_{k-1}!'$, with the maximum number of rows being $i_{k-1,max} = p_{k-1}!'$. The number of columns, as in the case of matrix $A_k$, can be arbitrarily large up to infinity.

Column numbers

| | | 1 | 2 | 3 | 4 | 5 | 6 | … |
|---|---|---|---|---|---|---|---|---|
| | 1 | | | | | | | |
| 1 | | 2 | 4 | 6 | 8 | 10 | 12 | … |
| 2 | | 3 | 5 | 7 | 9 | 11 | 13 | … |

*Matrix* $A_1$

Column numbers

| | | 1 | 2 | 3 | 4 | 5 | 6 | … |
|---|---|---|---|---|---|---|---|---|
| | 1 | | | | | | | |
| 1 | | 2 | 8 | 14 | 20 | 26 | 32 | … |
| 2 | | 3 | 9 | 15 | 21 | 27 | 33 | … |
| 3 | | 4 | 10 | 16 | 22 | 28 | 34 | … |
| 4 | | 5 | 11 | 17 | 23 | 29 | 35 | … |
| 5 | | 6 | 12 | 18 | 24 | 30 | 36 | … |
| 6 | | 7 | 13 | 19 | 25 | 31 | 37 | … |

*Matrix* $A_2$

Column numbers

row number

| | | 1 | 2 | 3 | 4 | 5 | 6 | … |
|---|---|---|---|---|---|---|---|---|
| | 1 | | | | | | | |
| 1 | | 2 | 32 | 62 | 92 | 122 | 152 | … |
| 2 | | 3 | 33 | 63 | 93 | 123 | 153 | … |
| 3 | | 4 | 34 | 64 | 94 | 124 | 154 | … |
| 4 | | 5 | 35 | 65 | 95 | 125 | 155 | … |
| 5 | | 6 | 36 | 66 | 96 | 126 | 156 | … |
| 6 | | 7 | 37 | 67 | 97 | 127 | 157 | … |
| 7 | | 8 | 38 | 68 | 98 | 128 | 158 | … |
| 8 | | 9 | 39 | 69 | 99 | 129 | 159 | … |
| 9 | | 10 | 40 | 70 | 100 | 130 | 160 | … |
| 10 | | 11 | 41 | 71 | 101 | 131 | 161 | … |
| 11 | | 12 | 42 | 72 | 102 | 132 | 162 | … |
| 12 | | 13 | 43 | 73 | 103 | 133 | 163 | … |
| 13 | | 14 | 44 | 74 | 104 | 134 | 164 | … |
| 14 | | 15 | 45 | 75 | 105 | 135 | 165 | … |
| 15 | | 16 | 46 | 76 | 106 | 136 | 166 | … |
| 16 | | 17 | 47 | 77 | 107 | 137 | 167 | … |
| 17 | | 18 | 48 | 78 | 108 | 138 | 168 | … |
| 18 | | 19 | 49 | 79 | 109 | 139 | 169 | … |
| 19 | | 20 | 50 | 80 | 110 | 140 | 170 | … |
| 20 | | 21 | 51 | 81 | 111 | 141 | 171 | … |
| 21 | | 22 | 52 | 82 | 112 | 142 | 172 | … |
| 22 | | 23 | 53 | 83 | 113 | 143 | 173 | … |
| 23 | | 24 | 54 | 84 | 114 | 144 | 174 | … |
| 24 | | 25 | 55 | 85 | 115 | 145 | 175 | … |
| 25 | | 26 | 56 | 86 | 116 | 146 | 176 | … |
| 26 | | 27 | 57 | 87 | 117 | 147 | 177 | … |
| 27 | | 28 | 58 | 88 | 118 | 148 | 178 | … |
| 28 | | 29 | 59 | 89 | 119 | 149 | 179 | … |
| 29 | | 30 | 60 | 90 | 120 | 150 | 180 | … |
| 30 | | 31 | 61 | 91 | 121 | 151 | 181 | … |

*Matrix* $A_3$

Column numbers

row number

| | | 1 | 2 | 3 | 4 | … |
|---|---|---|---|---|---|---|
| | 1 | | | | | |
| 1 | | 2 | 212 | 422 | 632 | … |
| 2 | | 3 | 213 | 423 | 633 | … |
| 3 | | 4 | 214 | 424 | 634 | … |
| 4 | | 5 | 215 | 425 | 635 | … |
| 5 | | 6 | 216 | 426 | 636 | … |
| 6 | | 7 | 217 | 427 | 637 | … |
| 7 | | 8 | 218 | 428 | 638 | … |
| 8 | | 9 | 219 | 429 | 639 | … |
| 9 | | 10 | 220 | 430 | 640 | … |
| 10 | | 11 | 221 | 431 | 641 | … |
| 11 | | 12 | 222 | 432 | 642 | … |
| 12 | | 13 | 223 | 433 | 643 | … |
| 13 | | 14 | 224 | 434 | 644 | … |
| 14 | | 15 | 225 | 435 | 645 | … |
| 15 | | 16 | 226 | 436 | 646 | … |
| 16 | | 17 | 227 | 437 | 647 | … |
| 17 | | 18 | 228 | 438 | 648 | … |
| 18 | | 19 | 229 | 439 | 649 | … |
| 19 | | 20 | 230 | 440 | 650 | … |
| 20 | | 21 | 231 | 441 | 651 | … |
| 21 | | 22 | 232 | 442 | 652 | … |
| 22 | | 23 | 233 | 443 | 653 | … |
| 23 | | 24 | 234 | 444 | 654 | … |
| 24 | | 25 | 235 | 445 | 655 | … |
| 25 | | 26 | 236 | 446 | 656 | … |
| 26 | | 27 | 237 | 447 | 657 | … |
| 27 | | 28 | 238 | 448 | 658 | … |
| 28 | | 29 | 239 | 449 | 659 | … |
| 29 | | 30 | 240 | 450 | 660 | … |
| 30 | | 31 | 241 | 451 | 661 | … |
| 31 | | 32 | 242 | 452 | 662 | … |
| 32 | | 33 | 243 | 453 | 663 | … |
| 33 | | 34 | 244 | 454 | 664 | … |
| 34 | | 35 | 245 | 455 | 665 | … |
| 35 | | 36 | 246 | 456 | 666 | … |
| 36 | | 37 | 247 | 457 | 667 | … |
| 37 | | 38 | 248 | 458 | 668 | … |
| 38 | | 39 | 249 | 459 | 669 | … |
| 39 | | 40 | 250 | 460 | 670 | … |
| 40 | | 41 | 251 | 461 | 671 | … |
| 41 | | 42 | 252 | 462 | 672 | … |
| … | | … | … | … | … | … |
| 207 | | 208 | 418 | 628 | 838 | … |
| 208 | | 209 | 419 | 629 | 839 | … |
| 209 | | 210 | 420 | 630 | 840 | … |
| 210 | | 211 | 421 | 631 | 841 | … |

*Matrix* $A_4$

**Fig.1. Matrices of prime numbers** [1]

[1] The number of rows of matrix $A_4$ exceeds the paper size, so here only a fragment is shown. Also, for this reason, other matrices are not shown in this figure.

In papers [2], and also in [1], to determine the properties of the developed matrices, we proposed a number of lemmas and theorems, as well as their detailed proofs. Therefore, for convenience, we present them without proof. In particular:

**Lemma 1.**

***Any integer $Z > 0$ occupies only one specific place in the matrix $A_k$. Moreover, the serial numbers $j_k$ (column number) and $i_k$(row number), in which a given number is located, are determined in a uniquely way:***

$$j_k = integer\ (\frac{Z-2}{p_k!'}) + 1$$

$$i_k = r_k + 1$$

***where $r_k$ - residue obtained by dividing the $(Z-2)$ on $p_k!'$, i.e.:***

$$(Z-2) = r_k(mod\ \ p_k!')$$

To determine whether an arbitrarily chosen number can occupy the same place in matrices $A_k$ and $A_{k+1}$, the following lemma was proved:

**Lemma 2.**

***If the randomly selected number takes a place located in the first column of the matrix $A_k$, then it occupies the same place in the matrix $A_{k+1}$, but in other cases it takes different places in matrices $A_k$ and $A_{k+1}$.***

Considering the set of numbers in an arbitrary row of matrix $A_k$, and to determine the distribution of these numbers in case of matrix $A_{k+1}$, we also proved the following lemma:

**Lemma 3.**

***The set of numbers located in one selected row of the matrix $A_k$, are redistributed in $p_{k+1}$ rows of the matrix $A_{k+1}$. At the same time, sequence of numbers contained in each of these $p_{k+1}$ rows, is an arithmetic progression with a constant $D_{k+1} = p_{k+1}!'$.***

This lemma automatically implies the proof of the following lemma:

**Lemma 4.**

***The set of numbers located in all the rows of the matrix $A_k$ are redistributed over $p_{k+1}!'$ rows of the matrix $A_{k+1}$.***

From Dirichlet's theorem for primes in an arithmetic progression it follows that if the first term and the difference of the progression are not mutually prime, then in this progression there will not be one prime number or there will be only one prime number. Moreover, this prime number is the first term of the arithmetic progression in question. It also follows from Dirichlet's theorem that if the first term and the difference of the arithmetic progression are mutually prime, then in this progression there is an infinite set of primes [3], [4] and [5].

In Fig.1 all of the rows in which there are only composite numbers (and possibly with one prime number), i.e. the first number and the difference in the arithmetic progression of which are not mutually prime, are, for clarity, repainted in a dark color and such rows will be called ordinary rows. The number of ordinary (repainted) rows will be **denoted** as $\boldsymbol{\beta_k}$,, where the index $k$ indicates belonging of these rows to the matrix $A_k$.

In this figure, all the rows in which there is an infinite number of prime numbers, i.e. the first number and the difference in the arithmetic progression of which are mutually prime, have been left uncolored for clarity. The total number of uncolored rows is **denoted** as $\boldsymbol{\alpha_k}$,, where the index $k$ also shows the belonging of these rows to the matrix $A_k$.

It follows from (1) that if $D_k > 2$, then prime twin numbers can not be located in a one particular uncolored row. Twins can only be in two adjacent uncolored rows, the first term and the difference in the arithmetic progression of which are mutually prime, and on the other hand the difference of the ordinal numbers of these two adjacent uncolored rows must be equal to 2. A pair of such uncolored rows is called *a pair of twin rows or twin-rows*. For two numbers located in different rows, but in one column of pairs of twin-rows, the equality $|a(k,i,j) - a(k,i \pm 2,j)| = 2$ is always valid. A pair of such prime numbers are twins. For example, from Fig. 1 it can be seen that prime twin numbers (5 and 7), (29 and 31), (41 and 43), and others are located in different rows, but in one column. The total number of pairs of twin rows is **denoted** by $\boldsymbol{\alpha_k^{twin}}$, where the index $k$ also shows the belonning of these rows to the matrix $A_k$.

If a row has an infinite set of prime numbers, but it is located from the same row at a distance greater than 2, then such row is called a *single row*. Therefore, in single rows there can not be any prime twin numbers. The total number of single rows is **denoted** as $\boldsymbol{\alpha_k^{single}}$, where the index $k$ also shows the belonging of these rows to the matrix $A_k$. Then the total number of uncolored rows is

$$\alpha_k = \alpha_k^{single} + 2\alpha_k^{twin} \tag{2}$$

If a total number of all rows to be **denoted** as $\boldsymbol{\Omega_k}$, then

$$\Omega_k = \alpha_k + \beta_k = \alpha_k^{single} + 2\alpha_k^{twin} + \beta_k \tag{3}$$

On the other hand,:

$$\Omega_k = i_{k,max} = D_k = p_k!' \tag{4}$$

The aim of this work is to determine the total number of primes-twins. Therefore, in the future we will be focusing on the pairs of twin-rows.

## 3. An infinite amount of prime twin numbers

As it was mentioned above, all twin prime numbers are mainly located in paired of twin-rows. If in this case at some point, for example, when considering $A_k$ matrix, all pairs of twin-rows disappear (i.e. they for some reason have become ordinary or single rows), then, it is obvious that they will not appear in any of the next matrices. In this case, it means that the number of twins should be limited.

Consider the number of rows of $A_k$ matrix in which prime numbers are being distributed.

**Theorem 1.** ***The number of unpainted rows of $A_k$ matrix, in which prime numbers are being distributed, is determined by the Euler's function $\varphi(D_k)$, where $D_k = p_k!'$- a constant in arithmetic progressions of numbers that are in each of these rows.***

**Proof of theorem 1.**

By definition the Euler's function $\varphi(D_k)$ determines the number of positive integers $l_k$, not exceeding $D_k$ and which are mutually prime[2] with $D_k$, i.e. $1 \leq l_k \leq D_k$ and $(D_k, l_k) \equiv 1$. On the other hand, from Dirichlet's theorem it follows that if the first term and the difference of arithmetic progression are mutually prime, then in this progression there are infinitely many prime numbers. Consequently, the number of rows of the matrix $A_k$, on which an infinite number of prime numbers is being distributed is equal to the value of Euler's function $\varphi(D_k)$:

$$\alpha_k = \varphi(D_k) = \varphi(p_k!') = (p_1 - 1)(p_2 - 1) \ldots (p_k - 1) = (p_k - 1)!' \qquad (5)$$

**Theorem 1 is proved.**

**Theorem 2.**

***Prime numbers, which are located in one unpainted row of the matrix $A_{k-1}$, are redistributed in $p_k - 1$ unpainted rows of matrix $A_k$.***

**Proof of theorem 2.**

From Theorem 1 it follows that in the case of the matrix $A_{k-1}$ prime numbers are redistributed according to its $\varphi(D_{k-1}) = \varphi(p_{k-1}!') = (p_{k-1} - 1)!'$ rows. In case of matrix $A_k$ redistribution of prime numbers occurs in $\varphi(D_k) = \varphi(p_k!') = (p_k - 1)!'$ rows. Consequently, if we consider the set of prime numbers that are in one of the unpainted rows of $A_{k-1}$ matrix, then their redistribution in the case of the matrix $A_k$ occurs over uncolored rows, the number of which is equal to

$$\frac{\varphi(D_k)}{\varphi(D_{k-1})} = \frac{(p_k - 1)!'}{(p_{k-1} - 1)!'} = p_k - 1 \qquad (6)$$

**Theorem 2 is proved.**

A number of corollaries follow from this theorem.

**Corollary 1 from Theorem 2.**

It follows from Lemma 3 that an infinite set of numbers in a chosen row (for example, in an uncolored) of the matrix $A_{k-1}$, are redistributed inside $p_k$ rows of the matrix $A_k$. On the other hand, from Theorem 2 we can state that an infinite set of prime numbers in the same uncolored row of matrix $A_{k-1}$, are redistributed inside $p_k - 1$ uncolored rows of matrix $A_k$.

---

[2] It follows from (1) and Fig. 1 that for our arithmetic progressions the condition $2 \leq l \leq D + 1$ is satisfied. In this case, if the number 1 is transferred to the matrix $A_k$, then all the members of the last row go to the 1st row and the penultimate row becomes the last row. Then the condition is $1 \leq l \leq D$ fulfilled automatically. Therefore, all the conclusions obtained for the conditions $1 \leq l \leq D$ are applicable and are used for our cases. We notice that in this case the first number of this last row will be equal $l = D = p_k!'$. Thus it is a row of composite numbers, which is a colored row.

*Consequently, the infinite set of some part of the composite numbers that were in an uncolored row of the matrix $A_{k-1}$ in the case of the matrix $A_k$ are separately redistributed as a separate ordinary (colored) row for which the first term and the difference of the arithmetic progression are not mutually prime.*

**Corollary 2 of Theorem 2.**

Now consider the general case. Suppose that there are $\alpha_{k-1}$ uncolored rows in a matrix $A_{k-1}$. Then, continuing the conclusions of Lemma 3, we find that the infinite set of numbers in these rows are redistributed inside $\alpha_{k-1}p_k$ rows of the matrix $A_k$. On the other hand, it also follows rom Theorem 2 that the infinite set of prime numbers in these $\alpha_{k-1}$ uncolored rows of the matrix $A_{k-1}$ are redistributed in $\alpha_{k-1}(p_k - 1)$ uncolored rows of matrix $A_k$.

*Consequently, the infinite set of some of the composite numbers that were in all of the uncolored $\alpha_{k-1}$ rows of the matrix $A_{k-1}$ , in the case of the matrix $A_k$ are separately redistributed in $\alpha_{k-1}$ separate ordinary (colored) rows, for each of which the first term and the difference of the arithmetic progression are not mutually prime.*

**Corollary 3 of Theorem 2.**

Now we determine exactly which composite numbers are separated from an infinite set of prime numbers and redistributed as separate colored rows in the matrix $A_k$.

As it was shown above, the first term and the difference of the progression of each of these colored rows are not mutually prime; from (1) we obtain that:

$$(l_k, D_k) \not\equiv 1 \text{ and } (l_k, p_k!') \not\equiv 1$$

This happens when the value of the first member of the arithmetic progression of the observed row is equal to <u>the value (or a multiple of that value) of one of the following primes $p_1,\ p_2,\ p_3, \dots ,\ p_k$.</u>

We will show this more clearly and in detail. For example, in the case of matrix $A_1$, only one prime number $p_1 = 2$ is being considered. Consequently, an infinite set of composite numbers, multiples of 2, which were in an infinite row of natural numbers, i.e. in the matrix $A_0$, in the case of the matrix $A_1$ are now separated from the prime numbers and arranged as a separate colored row (Fig. 1, $A_1$). Note that the first member of this row is 2, i.e. $l_1 = p_1 = 2$. Thus, during the formation of the matrix $A_1$, an infinite set of composite numbers that can be divided by 2 are fully identified and output in a separate row, which for clarity is recolored into a dark color (Fig. 1, $A_1$). Note that the prime number $p_1 = 2$ is also separated from the other prime numbers and being the only prime number is located in the given colored row along with other identified composite numbers.

In the case of matrix $A_2$, the next prime number is $p_2 = 3$. In this case, all of the numbers that in the case of the matrix $A_1$ were in its only uncolored row are redistributed along the three rows of the matrix $A_2$. In this case, one of these rows, in which the first term and the progression constant are not mutually prime, becomes colored. That is, for this row the following condition must be fulfilled:

$$(l_2, p_1 p_2) \not\equiv 1.$$

From this condition we obtain that the fulfillment of the requirement $(l_2, p_1) \not\equiv 1$ has already been considered in the case of matrix $A_1$. Therefore, for this row, the following condition must be fulfilled:

$$(l_2, p_2) \not\equiv 1.$$

This automatically implies that the first number of this colored row should be $l_2 = p_2 = 3$. The remaining numbers of this row must be equal to the multiple of the prime number $p_2 = 3$. Thus, it should be noted that during the formation of the $A_2$ matrix an infinite set of composite numbers that are divisible by 3 are fully revealed and are also output in a separate row, which for clarity is recolored into a dark color (Fig. 1, $A_2$). Note that in this case the number $p_2 = 3$ is also separated from other prime numbers and being the only prime number is located in the given colored row along with other identified composite numbers.

As a result, we get that in matrix $A_2$ an infinite set of composite numbers, which are divided by $p_1 = 2$ and $p_2 = 3$, are fully revealed and distributed in separate colored rows. Moreover, we also get that the prime numbers $p_1$ and $p_2$will not participate in the distribution of other prime numbers in the uncolored rows of the matrix $A_2$.

This pattern of distribution of prime and composite numbers is observed and occurs in all subsequent matrices. For example, continuing the above algorithm up to the matrix $A_k$, we get that an infinite set of composite numbers that are not mutually prime with the difference of the arithmetic progression $D_k = p_k!' = p_1 p_2 p_3 \dots p_k$, are distributed over the individual colored rows of the matrix $A_k$. It should be noted that these separated composite numbers can be conditionally divided into $k$ groups:

- The 1st group is an infinite set of composite numbers that are divisible by the number $p_1 = 2$,
- The 2nd group is an infinite set of composite numbers that are divisible by the number $p_2 = 3$,
- The 3rd group is an infinite set of composite numbers, which are divisible by the number $p_3 = 5$,

  …
- The $k$-th group is an infinite set of composite numbers, which are divisible by the number $p_k$.

Summarizing these groups, we can say that in the colored rows of the matrix $A_k$ an infinite set of those composite numbers that are not coprime with the numbers $p_1, p_2, p_3, \dots, p_k$ will be distributed.

It should be said that in this case, i.e. in the case of matrix $A_k$, the following prime numbers $p_1, p_2, p_3, \dots, p_k$ are also separated from other prime numbers and will not participate in the redistribution of an infinite set of other prime numbers in $(p_k - 1)!'$ rows of the matrix $A_k$. Each of them will be located in the corresponding colored rows of the matrix $A_k$, as the first number of this row and the corresponding arithmetic progression.

**Theorem 3.**

***The number of pairs of twin rows in the matrix $A_k$ increases monotonically with increasing order number k of the given matrix, and in each row of these pairs of twin rows there are an infinite number of primes.***

**Proof of Theorem 3.**

It follows from Lemma 3 that all the numbers in one uncolored row of the matrix $A_{k-1}$ (including in one of the rows of an arbitrary pair of twin rows) are redistributed over $p_k$ rows of the matrix $A_k$. Consequently, all the numbers in an arbitrary pair of twin rows of the matrix $A_{k-1}$ are redistributed in $p_k$ pairs of rows of the matrix $A_k$.

Consider the case, suppose there are $\alpha_{k-1}^{twin}$ pairs of twin rows in the matrix $A_{k-1}$.We choose one pair from them and will observe its first row.

It follows from Theorem 2 that all prime numbers in this row are redistributed in $p_k - 1$ rows of the matrix $A_k$. This means that, from the analyzed $p_k$ pairs of rows of the matrix $A_k$, one pair ceases to be a pair of twin rows.

If we consider the redistribution of the numbers in the second row of the given pair of twin rows of the matrix $A_{k-1}$, we similarly get that all the numbers in this row are redistributed along the other $p_k$ rows of the matrix $A_k$.
Moreover, the set of prime numbers will be redistributed along $p_k - 1$ rows of the $A_k$ matrix. This means that from the analyzed $p_k$ pairs of rows of the matrix $A_k$, one more pair ceases to be a pair of twin rows.

As a result, we obtain that all prime numbers in an arbitrary pair of twin rows of the matrix $A_{k-1}$, in the case of the matrix $A_k$, are redistributed along its $p_k - 2$ pairs of twin rows.

From the above, we find that from the analyzed $p_k$ pairs of rows of the matrix $A_k$, two pairs of rows cease to be a pair of twin rows. Thus, in this matrix there are 2 pairs, each of which consists of one colored (ordinary) and one uncolored row.

If we take into account that there are $\alpha_{k-1}^{twin}$ pairs of rows-twins in the matrix $A_{k-1}$, then the total number of pairs of twin rows in the matrix $A_k$ will be equal to:

$$\alpha_k^{twin} = \alpha_{k-1}^{twin} (p_k - 2) \qquad (7)$$

As it was shown above, there is only one pair of twin rows in the matrix $A_2$ (Fig. 1). With this in mind, it is easy to establish that $\alpha_k^{twin}$ the number of pairs of twin rows in matrix $A_k$ is determined by the following expression:

$$\alpha_k^{twin} = (p_2 - 2)(p_3 - 2) * \dots * (p_k - 2) = (p_k - 2)!', \text{where } k \geq 2 \qquad (8)$$

It follows from (8) that as the number $k$ of the matrix $A_k$ increases, the number of pairs of twin rows in it increases monotonically. On the other hand, the sequence of numbers in each row of $\alpha_k^{twin}$ pairs of twin rows is an arithmetic progression. Moreover, as it was shown above, the first term and the difference of each of these progressions are mutually prime. Thus, it follows from Dirichlet's theorem that every row of any pair of twin rows has an infinite quantity of prime numbers.

**Theorem 3 is proved.**

It follows from (8) that $\lim_{k\to\infty} \alpha_k^{twin} = \infty$, i.e. as k → ∞, the number of pairs of twin rows of the matrix $A_k$ tends towards infinity. If we assume that in each pair $\alpha_k^{twin}$ of twin rows there is at least one pair of prime twins, then we can definitely say that the number of twins is infinite. And this is possible when the prime numbers in each uncolored row, including rows of any pair of twin rows of the matrix $A_k$, are located closer and denser than in case of the matrix $A_{k-1}$.

In this case, from an infinite set of prime numbers in each of two rows of any pair of twin rows, at least 2 prime numbers will necessarily be in the same column and thus they will be twins. Yet, it really so?

To see this, we first analyze the following general case.

Suppose that in the matrix $A_{k-1}$ there are $\alpha_{k-1}^{single}$ single (uncolored) rows, $\beta_{k-1}$ ordinary (colored) rows and $\alpha_{k-1}^{twin}$ pairs of twin rows. Obviously, the total amount of rows is:

$$\Omega_{k-1} = \alpha_{k-1}^{single} + 2\alpha_{k-1}^{twin} + \beta_{k-1} = p_{k-1}!' \qquad (9)$$

Let's now define how many uncolored (single, as well as in pairs of twin rows) and colored (regular) rows will be in the matrix $A_k$.

It follows from Theorem 2 that all prime numbers in $\alpha_{k-1}^{single}$ single (non-colored) rows of the matrix $A_{k-1}$ are redistributed along $\alpha_{k-1}^{single}(p_k - 1)$ single rows of matrix $A_k$, i.e.:

$$\alpha_k^{single} = \alpha_{k-1}^{single}(p_k - 1) \qquad (10)$$

In addition, from this equality and Lemma 3 we also obtain that $\alpha_{k-1}^{single}$ rows become ordinary colored rows. The number of these newly formed colored rows is equal to:

$$\beta_k^s = \alpha_{k-1}^{single}. \qquad (11)$$

Here, the superscript *s* of the parameter $\beta_k^s$ shows that these colored (usual) rows of the matrix $A_k$ appear from single (uncolored) rows of matrix $A_{k-1}$.

From Lemma 3 we can also obtain that all numbers located in $\beta_{k-1}$ colored rows of the matrix $A_{k-1}$ are redistributed over $\beta_{k-1}\, p_k$ rows of the matrix $A_k$, i.e.:

$$\beta_k^b = \beta_{k-1}\, p_k \qquad (12)$$

Here, the superscript *b* of the parameter $\beta_k^b$ shows that these colored (ordinary) rows of the matrix $A_k$ appear from the ordinary (colored) rows of matrix $A_{k-1}$.

Let's now consider $\alpha_{k-1}^{twin}$ pairs of twin rows of the matrix $A_{k-1}$. From Theorem 3, and also from (7) and (8), we will find that the total number of pairs of twin rows in matrix $A_k$ is equal to:

$$\alpha_k^{twin} = \alpha_{k-1}^{twin}(p_k - 2) = (p_k - 2)!', \text{ where } k \geq 2.$$

At the same time, as shown above, some pairs of rows (total in $2\alpha_{k-1}^{twin}$) in the matrix $A_k$, cease to be pairs of twin rows: one row of each of these pairs of rows will appear to be single (uncolored) rows, and the other will appear to be ordinary (colored) rows. Hence, there appear new colored and uncolored rows in matrix $A_k$, and their total quantity is equal to:

$$\alpha_k^t = 2\alpha_{k-1}^{twin} \qquad (13)$$

$$\beta_k^t = 2\alpha_{k-1}^{twin} \quad (14)$$

Here, the superscript *t* of the parameters $\alpha_k^t$ and $\beta_k^t$ shows that these analyzed rows appear from pairs of twin rows of matrix $A_{k-1}$.

As a result of (8) - (14) we get that the total number of colored and uncolored rows in the matrix $A_k$ will be equal to:

$$\beta_k = \beta_k^s + \beta_k^b + \beta_k^t = \alpha_{k-1}^{single} + 2\alpha_{k-1}^{twin} + \beta_{k-1}\, p_k \quad (15)$$

$$\alpha_k = \alpha_k^{single} + \alpha_k^t + 2\alpha_k^{twin} = \alpha_{k-1}^{single}(p_k - 1) + 2\alpha_{k-1}^{twin}\, p_k - 2\alpha_{k-1}^{twin} \quad (16)$$

Using (15) and (16) we can find the total number of rows in the matrix $A_k$:

$$\Omega_k = \alpha_k + \beta_k = \alpha_{k-1}^{single} p_k + 2\alpha_{k-1}^{twin}\, p_k + \beta_{k-1}\, p_k \quad (17)$$

Comparing (9) and (17), we get:

$$\Omega_k = \Omega_{k-1} p_k = p_k!'$$

The last equality identically coincides with Lemma 4, which confirms the correctness of the above conclusions.

Now let's consider the number of prime and composite numbers located in uncolored rows of the matrix $A_k$ and not exceeding $x$ such that:

$$x \gg D_k = p_k!'. \quad (18)$$

Let $\pi_0$, $n_0$ and $m_0$ be a quantity of respectively prime, composite numbers, as well as all numbers contained in the arithmetic progression, consisting of a series of natural numbers[3] from *1* to *x*:

$$m_0 = \pi_0 + n_0.$$

Let $\pi_{ki}$, $n_{k,i}$ and $m_{k,i}$ be a quantity of respectively prime, composite numbers, as well as all numbers[4] not exceeding $x$ and located in one uncolored row with the serial number *i*:

$$m_{k,i} = \pi_{ki} + n_{k,i}$$

Then $\Pi_k$, $N_k$ and $M_k$ are total numbers of respectively prime, composite, as well as all available numbers not exceeding $x$ and located in **all** non-colored rows of the matrix $A_k$, are determined by simple summation:

$$\Pi_k = \sum_{i=1}^{\varphi(D_k)} \pi_{ki}$$

$$N_k = \sum_{i=1}^{\varphi(D_k)} n_{ki}$$

[3] An arithmetic progression consisting of a series of natural numbers from 1 to x is not part of the matrix $A_k$ considered (Fig. 1). Therefore, in order to describe elements of this progression, the index "0" is used.

[4] Numerous scientific works are devoted to the problem of prime numbers [4]. For example, in 1935 Siegel obtained that for any fixed A> 1 for $1 \le D_k \le \log^A x$, for an arithmetic progression with parameters $(D, l) \equiv 1$ и $1 \le l \le D$, it is true that:

$$\pi(x, D, l) = \frac{Li(x)}{\varphi(D)} + O\left(x \exp\left(-c\sqrt{\ln x}\right)\right).$$

where D is the difference of the progression, *l* is the first term of the progression. This expression in the scientific literature is known as the Siegel-Walphlisch formula. In this formula, $\varphi(D)$ is the Euler function, $A$ is a positive constant, $c = c(A) > 0$ is a constant, where $\mathrm{Li}(x) = \int_2^x \frac{du}{\ln u}$. At the present time, many functions are known regarding the quantity of prime numbers not exceeding x in an arithmetic progression. Moreover, if the generalized Riemann hypothesis is considered correct, then, as known, the following formula holds:

$$\pi(x, D, l) = \frac{Li(x)}{\varphi(D)} + O\left(\sqrt{x}\,\ln x\right)$$

It should be said that the last function is the closest to the real values of $\pi(x, D, l)$. Nevertheless, many problems related to even more exact approximations remain unresolved [4], [5].

$$M_k = \Pi_k + N_k = \sum_{i=1}^{\varphi(D_k)}(\pi_{ki} + n_{ki})$$

Here the summation is carried out over uncolored rows.

Let's introduce several new notation terms. Let $\pi_{k,av}$ , $n_{k,av}$ and $m_{k,av}$ be average numbers of prime, composite and all numbers not exceeding $x$ and located in one <u>uncolored</u> row of matrix $A_k$. Then,

$$\pi_{k,av} = \frac{\Pi_k}{\varphi(D_k)} = \frac{\sum_{i=1}^{\varphi(D_k)} \pi_{k,i}}{\varphi(D_k)}$$

$$n_{k,av} = \frac{N_k}{\varphi(D_k)} = \frac{\sum_{i=1}^{\varphi(D_k)} n_{k,i}}{\varphi(D_k)}$$

$$m_{k,av} = \frac{\Pi_k + N_k}{\varphi(D_k)} = \frac{\sum_{i=1}^{\varphi(D_k)}(\pi_{k,i} + n_{ki})}{\varphi(D_k)} = \pi_{k,av} + n_{k,av} \qquad (19)$$

On the other hand, it is obvious that the total number of all numbers (prime and composite) not exceeding $x$ and located in any row of the matrix $A_k$ are equal to each other and defined by the relation:

$$m_{k,av} = m_{k,i} = \frac{x}{D_k} = \frac{x}{p_k!'}$$

Here and later on we shall consider the cases when condition (18) is satisfied, i.e. $x \gg D_k = p_k!'$.

Let's determine the average density of prime numbers not exceeding $x$ and located in one <u>uncolored</u> row of the matrix $A_k$:

$$\rho_{k,av} = \frac{\pi_{k,av}}{m_{k,av}} = \frac{\pi_{k,av} p_k!'}{x} \qquad (20)$$

**Theorem 4.**

***The density of prime numbers in uncolored rows of the matrix $A_k$ grows monotonically with the increase in the sequence number k of the observed matrix. And the density of prime numbers in each uncolored row in limit tends towards the limit value:***

$$\lim_{k\to\infty} \rho_k = 1$$

**Proof of theorem 4.**

It follows from Theorem 1 that the prime numbers that are distributed along $\varphi(D_{k-1})$ rows of the matrix $A_{k-1}$, in the case of matrix $A_k$, are being redistributed along its $\varphi(D_k)$ rows. On the other hand, from the corollary to Theorem 2, we obtain that this redistribution occurs without any participation of a prime number $p_k$. Therefore, the following equalities hold:

$$\pi_1\varphi(D_1) = \pi_0\varphi(D_0) - 1$$
$$\pi_{2,av}\varphi(D_2) = \pi_{1,av}\varphi(D_1) - 1$$
$$\pi_{3,av}\varphi(D_3) = \pi_{2,av}\varphi(D_2) - 1$$

……………………..………….……..

$$\pi_{k,av}\varphi(D_k) = \pi_{k-1.av}\varphi(D_{k-1}) - 1$$

From the last equality we get that:

$$\pi_{k,av} = \frac{\pi_{k-1,av}\varphi(D_{k-1})}{\varphi(D_k)} - \varepsilon_1 = \frac{\pi_{k-1,av}}{p_k - 1} - \varepsilon_1, \qquad (21)$$

where the small value $\varepsilon_1 = \frac{1}{\varphi(D_k)}$. From (20) and (21) we obtain that:

$$\rho_{k,av} = \rho_{k-1,av} \frac{p_k}{(p_k - 1)} - \varepsilon_2, \qquad (22)$$

where the small value $\varepsilon_2 = \frac{\varepsilon_1}{m_{k,av}} = \frac{1}{\varphi(D_k) m_{k,av}}$.

From (22) it follows that:

- prime numbers in any uncolored row of matrix $A_k$ are more densely distributed than in case of the previous matrix $A_{k-1}$,
- with an increase in order number $k$ of the matrix $A_k$, the density of prime numbers in any of its uncolored rows increases continuously.

On the other hand, as $k \to \infty$, it becomes obvious that:

$$\lim_{k \to \infty} \frac{p_k}{(p_k - 1)} = 1$$

Consequently, taking this into account, we obtain from (22) that:

$$\lim_{k \to \infty} \frac{\rho_{k,av}}{\rho_{k-1.av}} = 1, \qquad (23)$$

i.e. as $k \to \infty$ in the limit, we get that $\lim_{k \to \infty} \rho_{k,av} = \lim_{k \to \infty} \rho_{k-1,av}$.

From (19), (20) and (23) we get that:

$$\lim_{k \to \infty} \frac{\rho_{k,av}}{\rho_{k-1,av}} = \lim_{k \to \infty} \frac{\frac{\pi_{k,av}}{m_k}}{\frac{\pi_{k-1,av}}{m_{k-1}}} = \lim_{k \to \infty} \frac{\pi_{k,av}}{\pi_{k-1,av}} * \frac{\pi_{k-1,av} + n_{k-1,av}}{\pi_{k,av} + n_{k,av}} = 1 \qquad (24)$$

Let us analyze equality (24). If as $k \to \infty$ the following equality holds:

$$\lim_{k \to \infty} n_{k,av} = \lim_{k \to \infty} n_{k-1,av} = 0, \qquad (25)$$

then (24) is automatically satisfied.

Consequently, it follows from (19), (20) and (25) that

$$\lim_{k \to \infty} \rho_{k,av} = \lim_{k \to \infty} \frac{\pi_{k,av}}{m_{k,av}} = \lim_{k \to \infty} \frac{\pi_{k,av}}{\pi_{k,av} + n_{k,av}} = \lim_{k \to \infty} \frac{\pi_{k,av}}{\pi_{k,av}} = 1 \qquad (26)$$

Equation (26) means that as the order number of the matrix $A_k$ increases, the density of prime numbers in its uncolored rows grows monotonically and tends to the limiting value.

**Theorem 4 is proved.**

It should be said that the fulfillment of Theorem 4, in particular (26), is possible when prime numbers in any uncolored row of matrix $A_k$ are more densely distributed than in case of uncolored rows of previous matrices. This means that in each pair of twin rows of the matrix $A_k$ there will undoubtedly be twins.

Now let’s go back to the main problem of this work.

**Theorem 5.** ***The hypothesis of an infinite quantity of prime twins is true.***

**Proof of theorem 5.**

Consider the matrix $A_k$ for a truly large $k \to \infty$. It follows from Dirichlet's theorem that every row of any pair of twin rows of $A_k$ contains an infinite set of primes.

On the other hand, it follows from Theorem 3 that as the order number of the matrix $A_k$ grows, the number of pairs of twin rows in it grows monotonically, such that

$$\lim_{k\to\infty} \alpha_k^{twin} = \lim_{k\to\infty} (p_k - 2)!' = \infty$$

It also follows from Theorem 4 and (26) that as the order number of the matrix $A_k$ grows, the density of prime numbers in each uncolored row grows monotonically. This means that in each pair of twin rows there are twin numbers. This, in turn, implies that the number of twins is infinite.

**Theorem 5 is proved.**

## 4. Conclusion.

In this paper, in order to study the properties of prime twin numbers, we first introduce the notion of a matrix of prime numbers. Then we consider the density of prime numbers in separate rows of a matrix to be considered. In this paper we prove a number of lemmas and theorems that together with the Dirichlet and Euler theorems make it possible to prove the infinity of prime twins.

**P.S. The authors would be sincerely grateful to those who will send their comments and suggestions to: baibekovsn@mail.ru**

**AUTHOR'S BIOGRAPHY**

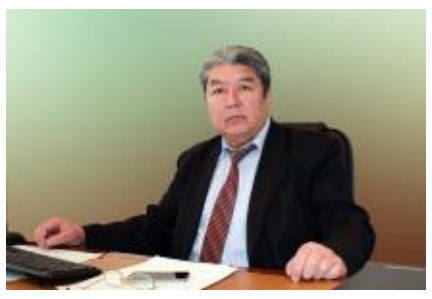

Seidikassym N. Baibekov
Nationality - Republic of Kazakhstan
Education – Physical faculty of Leningrad State University
Academic degree - Doctor of technical science
Academic title- Professor
Author of more than 130 scientific articles and about 20 monographs, textbooks and manuals
Academic interests - mathematics (number theory), physics, IT-technology and others.